\documentclass[twoside,fleqn, 12pt]{article}
\usepackage{amsmath,amssymb,amsfonts, epsfig}
\newtheorem{theorem}{Theorem} 
\newtheorem{corollary}{Corollary} \newtheorem{proposition}{Proposition}
 
\newtheorem{remark}{Remark}
\numberwithin{equation}{section} \numberwithin{lemma}{section}
\numberwithin{theorem}{section} \numberwithin{corollary}{section}
\numberwithin{proposition}{section}
\numberwithin{definition}{section} \numberwithin{example}{section}
\numberwithin{remark}{section}
\def\TITLE       {\textbf{About some split central simple algebras}}% <--- Title of Talk\\
\par
\bigskip
\bigskip

\def\FIRSTNAMEI  {Diana} % <--- First Author's First Name
\def\LASTNAMEI {Savin\\}% <--- First Author's Last Name\\
\begin{document}

\def\ABSTRACT {
\begin{center}
\textbf{Abstract}
\end{center}
\begin{small}
{In this paper we study certain quaternion algebras and symbol algebras which split.}
\end{small}
 }% }%<-- Optional
\def\CLASSIFICATION{MSC (2010): 11A41, 11R04, 11R18, 11R37, 11R52, 11S15, 11F85}% <--- Optional

\def\KEYWORDS{KEYWORDS: quaternion algebras; symbol algebras; cyclotomic fields; Kummer fields, p-adic fields}% <-- Optional
%\frontmatter++++++++++++++++++++++++++++++++++++++++++++++++++++%
\begin{center}
\par\TITLE
\par\FIRSTNAMEI \ \LASTNAMEI
\par%\medskip
%Elena - aici am pus universitatea
\par

\end{center}
\par\
\par\ABSTRACT{}
\par\
\par\CLASSIFICATION{}
\par\KEYWORDS{}
\par\

%++++++++++++++++++++++++++++++++++++++++++++++++++++++++++++++++++++++++++
%\frontmater

\section{Introduction}
\label{Introduction}

Let $K$ be a field such that all vector spaces over $K$ are finite
dimensional. Let $A$ be a simple $K$-algebra and $Z\left(A\right)$ be the centre of $A.$ We recall that the $K$- algebra $A$ is called central simple if $Z\left(A\right)=K.$\\
Let $K$ be a field with $charK\neq 2.$ Let $\Bbb{H}_{K}\left( \alpha
,\beta \right) $ be the generalized quaternion algebra with basis $\{1,e_{1},e_{2},e_{3}\}\,$%
\thinspace and the multiplication given by

{\footnotesize
\[
\begin{tabular}{l|llll}
$\cdot $ & $1$ & $\,\,\,e_{1}$ & $\,\,e_{2}$ & $\,\,\,\,e_{3}$ \\
\hline
$1$ & $1$ & $e_{1}$ & $e_{2}$ & $\,\,\,e_{3}$ \\
$e_{1}$ & $e_{1}$ & $\alpha $ & $e_{3}$ & $\alpha e_{2}$ \\
$e_{2}$ & $e_{2}$ & $-e_{3}$ & $\beta $ & $-\,\,\,\beta e_{1}$ \\
$e_{3}$ & $e_{3}$ & $-\alpha e_{2}$ & $\beta e_{1}$ & $-\alpha \beta $%
\end{tabular}
.
\]
}
A natural generalization of the quaternion algebras are the symbol algebras.
Let $n$ be an arbitrary positive integer, let \ $K$ be a field
whose $char(K)$ does not divide $n$ and contains a primitive $n$-th root of
unity. Denote $K^{\ast }=K\backslash \{0\},$ $a,b$ $\in K^{\ast }$ and let $%
S$ be the algebra over $K$ generated by elements $x$ and $y,$ where%
\begin{equation*}
x^{n}=a,y^{n}=b,yx=\xi xy.
\end{equation*}
This algebra is called a \textit{symbol algebra } and it is denoted by $\left( \frac{a,~b}{%
K,\xi }\right) .$ J. Milnor, in [19], calls it the symbol algebra
because of its connection with the $K$-theory and with the Steinberg symbol.
For $n=2,$ we obtain the quaternion algebra.\\
Quaternion algebras and symbol algebras are central simple algebras.\\
Quaternion algebras and symbol algebras have been studied from several points of view: from the theory of associative algebras ([20], [12], [6], [7], [9], [16], [23]), from number theory ([18], [12], [21], [15]), analysis and mecanics ([14]).\\
In this paper, we will determine certain split quaternion algebras and split symbol algebras, using some results of number theory (ramification theory in algebraic number fields, class field theory).\\
Let $K\subset L$ be a fields extension and let $A$ be a central simple algebra over the field $K.$ We recall that $A$ is called split by $L$ and $L$ is called a splitting field for $A$ if $A\otimes_{K}L$ is a matrix algebra over $L.$\\
In [12] appear the following criterions to decide if a quaternion algebra or a symbol algebra is split.
\begin{proposition}\label{1.1.} 
\textit{The quaternion algebra }$\Bbb{H}%
_{K}\left( \alpha ,\beta \right) $\textit{\ is split if and only if
the
conic }$C\left( \alpha ,\beta \right):$ $\alpha
x^{2}+\beta y^{2}=z^{2}$ \textit{\ has a rational point over }%
$K($\textit{i.e. if there are }$x_{0},y_{0},z_{0}\in K$\textit{\ such that }$%
\alpha x_{0}^{2}+\beta y_{0}^{2}=z_{0}^{2}).$
\end{proposition}

\begin{theorem}\label{1.}
 \textit{Let }$K$\textit{\ be a field such that }$%
\zeta \in K,\,\,\zeta ^{n}=1,\zeta $\textit{\ is a primitive root, and let }$%
\alpha ,\beta \in K^{*}.$\textit{\ Then the following statements are
equivalent:}\\
\textit{i) The cyclic algebra }$A=\left( \frac{\alpha ,\beta }{K,\zeta }%
\right) $\textit{\ is split.}\\
\textit{ii) The element }$\beta $\textit{\ is a norm from the
extension }$K \subseteq K(\sqrt[n]{\alpha}). $
\end{theorem}

\begin{remark}\label{1}([16]) \textit{Let} $K$ \textit{be an algebraic numbers field such that} $\left[K:\mathbb{Q}\right]$ \textit{is odd and} $\alpha ,\beta$$\in$$\mathbb{Q^{*}}.$ \textit{Then, the quaternion algebra} $\Bbb{H}%
_{K}\left( \alpha ,\beta \right)$ \textit{splits if and only if} \textit{the quaternion algebra} $\Bbb{H}%
_{\mathbb{Q}}\left( \alpha ,\beta \right)$ \textit{splits}.
\end{remark}

If in Proposition 1.1 we have $K=\mathbb{Q},$ to decide if 
the conic $C\left( \alpha ,\beta \right):$ $\alpha
x^{2}+\beta y^{2}=z^{2}$ has a rational point, we will use Minkovski-Hasse theorem.\\
\vspace{2mm}\\
\textbf{Minkovski-Hasse Theorem.} ([4]) \textit{A quadratic form with rational coefficients represents zero in the field of rational numbers if and only if it represents zero in the field of real numbers and in all fields of} $p$-\textit{adic numbers (for all primes} $p$).\\
For a quadratic form in three variables, Minkovski-Hasse Theorem can be reformulate as the following:\\
\textit{The form with rational coefficients} $\alpha x^{2}+\beta y^{2}- z^{2}$ \textit{with nonzero rational coefficients} $\alpha$ \textit{and} $\beta$ \textit{represents zero in the field of rational numbers if and only if for all primes} $p$ (\textit{including} $p=\infty$), we have

$$\left(\frac{\alpha,\beta}{p}\right)=1,$$
\textit{where} $\left(\frac{\alpha,\beta}{p}\right)$ \textit{is the Hilbert symbol in the} $p$-\textit{adic field} $\mathbb{Q}_{p}.$
\begin{corollary}\label{1.} ([4])\\
i) If $p$ is not equal with $2$ or $\infty$ and $p$ does not enter into the factorizations of $\alpha$ and $\beta$ into prime powers (which means that $\alpha$ and $\beta$ are $p$-adic units), then the form $\alpha x^{2}+\beta y^{2}- z^{2}$ represents zero in the $p$-adic fields $\mathbb{Q}_{p}$ and thus for all such $p$ the Hilbert symbol $\left(\frac{\alpha,\beta}{p}\right)=1.$\\
ii) $\left(\frac{\alpha,\beta}{\infty}\right)=1,$ if $\alpha>0$ or $\beta>0,$\\
$\left(\frac{\alpha,\beta}{\infty}\right)=-1,$ if $\alpha<0$ and $\beta<0,$\\
where $\left(\frac{\alpha,\beta}{\infty}\right)$ \textit{is the Hilbert symbol in the field} $\mathbb{R}.$
\end{corollary}
\begin{corollary}\label{2} ([4])\\
The product of the Hilbert symbols in the $p$-adic fields satisfy
$$\prod_{p}\left(\frac{\alpha,\beta}{p}\right)=1,$$
\textit{where} $p$-\textit{runs through all prime numbers and the symbol} $\infty.$
\end{corollary}
Now we recall a result about primes of the form $q=x^{2}+ny^{2}$ which we will use for study the quaternion algebras.
\begin{theorem}\label{theorem 1.2.} ([5], [11], [22], [3]) For an odd prime positive integer $q$,  the following statements are true:\\
i) $q=x^{2}+3y^{2}$ for some $x,y$ $\in$ $\mathbb{Z}$ if and only if $q\equiv1$(mod $3$) or $q=3;$\\
ii) $q=x^{2}+5y^{2}$ for some $x,y$ $\in$ $\mathbb{Z}$ if and only if $q\equiv 1;9$(mod $20$) or $q=5;$\\
iii) $q=x^{2}+6y^{2}$ for some $x,y$ $\in$ $\mathbb{Z}$ if and only if $q\equiv1;7$(mod $24$);\\
iv) $q=x^{2}+7y^{2}$ for some $x,y$ $\in$ $\mathbb{Z}$ if and only if $q\equiv1;2;4$(mod $7$) or $q=7$;\\
v) $q=x^{2}+10y^{2}$ for some $x,y$ $\in$ $\mathbb{Z}$ if and only if $q\equiv1;9;11;19$(mod $40$);\\
vi) $q=x^{2}+13y^{2}$ for some $x,y$ $\in$ $\mathbb{Z}$ if and only if $q\equiv1;9;17;25;29;49$(mod $52$) or $q=13$;\\
vii) $q=x^{2}+14y^{2}$ for some $x,y$ $\in$ $\mathbb{Z}$ if and only if $q\equiv1;9;15;23;25;39$(mod $56$);\\
viii) $q=x^{2}+15y^{2}$ for some $x,y$ $\in$ $\mathbb{Z}$ if and only if $q\equiv1;19;31;49$(mod $60$);\\
ix) $q=x^{2}+ 21y^{2}$ for some $x,y$ $\in$ $\mathbb{Z}$ if and only if $q\equiv1;25;37$(mod $84$);\\
x) $q=x^{2}+ 22y^{2}$ for some $x,y$ $\in$ $\mathbb{Z}$ if and only if $q\equiv1;9;15;23;25;31;47;49;71;81$(mod $88$);\\
xi) $q=x^{2}+ 30y^{2}$ for some $x,y$ $\in$ $\mathbb{Z}$ if and only if $q\equiv1;31;49;79$(mod $120$).
\end{theorem}
We recall some properties of cyclotomic fields and Kummer fields which will be necessary in our proofs.
\begin{proposition}\label{p1.2.} ([2]) Let $q$ be an odd positive prime integer and $\xi$ be a primitive root of order $q$ of the unity. Then the ring $\mathbb{Z}\left[\xi\right]$ is a principal domain for $q$$\in$$\left\{3,5,7,11,13,17,19\right\}.$
\end{proposition}
\begin{proposition}\label{p1.3.} ([13]) \textit{Let }$l$\textit{\ be a natural number, }$%
l\geq 3$\textit{\ and }$\zeta $\textit{\ be a primitive root of the unity
of\thinspace \thinspace }$l$-\textit{order. If }$p$\textit{\ is a prime
natural number, }$l$\textit{\ is not divisible with }$p$\textit{\ and }$f^{'}$%
\textit{\ is the smallest positive integer such that }$p^{f^{'}}\equiv 1$\textit{%
\ mod }$l$\textit{, then we have} 
\[
p\Bbb{Z}[\zeta ]=P_{1}P_{2}....P_{r}, 
\]
\textit{where }$r=\frac{%
\begin{array}{c}
\\ 
\varphi \left( l\right)
\end{array}
}{%
\begin{array}{c}
f^{'}
\end{array}
},\varphi $\textit{\ is the Euler's function and }$P_{j},\,j=1,...,r$\textit{%
\ are different prime ideals in the ring }$\Bbb{Z}[\zeta ].$%
\[
\]
\end{proposition}
\begin{proposition}\label{p1.4.} ([13])
\textit{Let }$\xi $\textit{\ be a primitive
root of the unity of }$q-$\textit{order, where }$q$\textit{\ is a prime
natural number and let} $A$ \textit{be the ring of integers of the Kummer
field} $Q(\xi, \sqrt[q]{\mu})$ . \textit{A prime ideal} $P$\textit{\ in
the ring} $\Bbb{Z}[\xi ]$\textit{\ is in }$A$\textit{\ in one of the
situations:}

\textit{i) It is equal with the }$q-$\textit{power of a prime ideal from} $A,
$ \textit{if the} $q-$\textit{power character }$\left( \frac{\mu }{P}\right)
_{q}=0;$

\textit{ii) It is a prime ideal in} $A$, \textit{if} $\left( \frac{\mu }{P}%
\right) _{q}=$\textit{\ a root of order }$q$\textit{\ of unity, different
from }$1$. \\
\textit{iii) It decomposes in $q$ different prime ideals from} $A$, \textit{%
if} $\left( \frac{\mu }{P}\right) _{q}=1.$
\end{proposition}

\begin{theorem}\label{theorem 1.3}([1],[15])) Let $K$ be an algebraic
number field, $v$ be a prime of $K$ and $K\subseteq L$ be 
a Galois extension. Let $w$ be a prime of $L$ 
lying above $v$ such that $K_{v}\subseteq L_{w}$
is a unramified extension of $K_{v}$ of (residual) degree $f.$ 
Let $b=\pi _{v}^{m}\cdot u_{v}\in K_{v}^{*},$ where $
\pi_{v}$ denote a prime element in
 $ K_{v}$ and $u_{v}$ a unit in the ring of integers $\mathcal{O}_{v},m\in \Bbb{Z}%
.$ Then $b\in N_{L_{w}\,/\,K_{v}}\left( L_{w}^{*}\right) $
 if and only if $f\,\mid \,m.$ In particular, every unit of $%
\mathcal{O}_{v}$ is the norm of a unit in $L_{w}.$
\end{theorem}

\section{\label{results} Main Results}

In the paper [21] we proved the following Propositions:

\begin{proposition}\label{p2.1.}
For $\alpha =-1,\beta =q,$ where $q$ is %
 a prime number, $q\equiv$$3$ (mod $4$), $K=\Bbb{Q},$ the algebra $\Bbb{H}_{\Bbb{Q}%
}\left( -1,q\right) $ is a division algebra.
\end{proposition}

\begin{proposition}\label{p2.2.}
\textit{If }$K=$\textit{\ }$\Bbb{Q}\left( \sqrt{3}%
\right) ,$\textit{\ then the quaternion algebra }$\Bbb{H}_{K}\left(
-1,q\right) ,$\textit{\ where } $q$ is %
 a prime number, $q\equiv 1$ (\textit{\ mod
}$3$), \textit{\ is a split algebra.}
\end{proposition}
In [16] appear the following results:
\begin{proposition}\label{p2.3.}
If $q$ is an odd prime positive integer, then:\\
i) the algebra $\Bbb{H}_{\Bbb{Q}%
}\left( -1,q\right) $ is a split algebra if and only if $q\equiv 1$ (mod $4$).\\
ii) the algebra $\Bbb{H}_{\Bbb{Q}%
}\left( -2,q\right) $ is a split algebra if and only if $q\equiv 1$ or $3$ (mod $8$).
\end{proposition}
In what follows we will give a sufficient condition such that the algebras $\Bbb{H}_{\Bbb{Q}%
}\left(\alpha,q\right),$ where $\alpha$$\in$$\left\{-3,-5,-6,-7,-10,-13,-14,-15,-21,-22,-30\right\}$ are split algebras.
\begin{proposition}\label{p2.4}Let $q$ be a prime positive integer. The following statements holds true:\\
i)if $q\equiv1$(mod $3$) or $q=3,$ then the algebra $\Bbb{H}_{\Bbb{Q}%
}\left( -3,q\right) $ is a split algebra;\\
ii)if $q\equiv1;9$(mod $20$) or $q=5,$ then the algebra $\Bbb{H}_{\Bbb{Q}%
}\left(-5,q\right) $ is a split algebra;\\
iii) if $p\equiv1;7$(mod $24$), then the algebra $\Bbb{H}_{\Bbb{Q}%
}\left(-6,q\right) $ is a split algebra;\\
iv) if $p\equiv1;2;4$(mod $7$) or $q=7$, then the algebra $\Bbb{H}_{\Bbb{Q}%
}\left(-7,q\right) $ is a split algebra;\\
v) if $q\equiv 1; 9; 11; 19$(mod $40$), then the algebra $\Bbb{H}_{\Bbb{Q}%
}\left(-10,q\right) $ is a split algebra;\\
vi) if $q\equiv 1; 9; 17; 25; 29; 49$(mod $52$), then the algebra $\Bbb{H}_{\Bbb{Q}%
}\left(-13,q\right) $ is a split algebra;\\
vii) if $q\equiv 1; 9; 15; 23; 25; 39$(mod $56$), then the algebra $\Bbb{H}_{\Bbb{Q}%
}\left(-14,q\right) $ is a split algebra;\\
viii) if $q\equiv 1; 19; 31; 49$(mod $60$), then the algebra $\Bbb{H}_{\Bbb{Q}%
}\left(-15,q\right) $ is a split algebra;\\
ix) if $q\equiv 1; 25; 37$(mod $84$), then the algebra $\Bbb{H}_{\Bbb{Q}%
}\left(-21,q\right) $ is a split algebra;\\
x) if $q\equiv 1; 9; 15; 23; 25; 31; 47; 49; 71; 81$(mod $88$), then the algebra $\Bbb{H}_{\Bbb{Q}%
}\left(-22,q\right) $ is a split algebra;\\
xi) if $q\equiv 1; 31; 49; 79$(mod $120$), then the algebra $\Bbb{H}_{\Bbb{Q}%
}\left(-30,q\right) $ is a split algebra.
\end{proposition}
\textbf{Proof.} The proof is immediately using Proposition 1.1 and Theorem 1.2.\\
\vspace{2mm}\\
We asked ourselves if the converse statements in Proposition 2.4 are true. We obtained that these are true for the algebras
$\Bbb{H}_{\Bbb{Q}%
}\left( -3,q\right),$ $\Bbb{H}_{\Bbb{Q}%
}\left( -5,q\right), $ $\Bbb{H}_{\Bbb{Q}%
}\left( -7,q\right), $ $\Bbb{H}_{\Bbb{Q}%
}\left( -13,q\right). $

\begin{proposition}\label{p2.5.}
Let $q$ be a prime positive integer. The following statements holds true:\\
i) the algebra $\Bbb{H}_{\Bbb{Q}%
}\left( -3,q\right) $ is a split algebra if and only if $q\equiv1$(mod $3$) or $q=3$;\\
ii) the algebra $\Bbb{H}_{\Bbb{Q}%
}\left(-5,q\right) $ is a split algebra if and only if $q\equiv1,9$(mod $20$) or $q=5$;\\
iii)the algebra $\Bbb{H}_{\Bbb{Q}%
}\left(-7,q\right) $ is a split algebra if and only if if $p\equiv1,2,4$(mod $7$) or $q=7$;\\
iv) the algebra $\Bbb{H}_{\Bbb{Q}%
}\left(-13,q\right) $ is a split algebra if and only if $q\equiv 1; 9; 17; 25; 29; 49$\\
(mod $52$).
\end{proposition}
\textbf{Proof.} For i), ii), iii), iv) we proved the implication $"\Leftarrow"$ in Proposition 2.4\\
We prove the implication $"\Rightarrow"$ only for i) and iv) (proofs for ii) and iii) are similar).\\
i) $"\Rightarrow"$ If the algebra 
$\Bbb{H}_{\Bbb{Q}%
}\left( -3,q\right) $ is a split algebra, applying Proposition 1.1 it results that the form $-3x^{2}+qy^{2}-z^{2}$ represents zero in the field of rational numbers. According to Minkovski-Hasse Theorem, this is equivalent with the form $-3x^{2}+qy^{2}-z^{2}$ represents zero in the field of real numbers and in all fields of $p$-adic numbers. The last statement is equivalent to the 
Hilbert symbol $\left(\frac{-3,q}{p}\right)=1,$ for all primes $p$ (including $p=\infty$).\\
According to Corollary 1.1 the Hilbert symbols $\left(\frac{-3,q}{\infty}\right)=1$ and $\left(\frac{-3,q}{p}\right)=1,$ for all primes $p\neq 2,3,q.$\\
\textbf{Case 1}: if $q=3.$\\
Similarly with i) The algebra 
$\Bbb{H}_{\Bbb{Q}%
}\left( -3,3\right) $ is a split algebra
if and only if the form $-3x^{2}+3y^{2}-z^{2}$ represents zero in the field of rational numbers.This is true, a solution is $\left(x_{0}, y_{0}, z_{0}\right)=\left(1, 1, 0\right).$\\
\textbf{Case 2}: if $q\neq 3.$\\
We detemine the values of $q$ for which the Hilbert symbols $\left(\frac{-3,q}{q}\right)=1,$ $\left(\frac{-3,q}{3}\right)=1$ and
$\left(\frac{-3,q}{2}\right)=1$
Using the properties of the Hilbert symbol we have:
$$\left(\frac{-3,q}{q}\right)=\left(\frac{-1,q}{q}\right)\cdot \left(\frac{3,q}{q}\right)=\left(\frac{-1}{q}\right)\cdot \left(\frac{3}{q}\right)=\left(-1\right)^{\frac{q-1}{2}}\cdot \left(\frac{3}{q}\right).$$
Applying Reciprocity law we obtain rapidly that the Hilbert symbol $\left(\frac{-3,q}{q}\right)=1$ if and only if $q\equiv 1$ (mod $3$).\\
$\left(\frac{-3,q}{2}\right)=\left(-1\right)^{\frac{-3-1}{2}\cdot\frac{q-1}{2}}=1.$\\
Using Corollary 1.2 or direct calculations we obtain that
 $\left(\frac{-3,q}{3}\right)=1$ if and only if $q\equiv 1$ (mod $3$).\\
From the previously proved and from Proposition 2.4, it results that the algebra $\Bbb{H}_{\Bbb{Q}%
}\left( -3,q\right) $ is a split algebra if and only if $q\equiv 1$ (mod $3$).\\
iv) $"\Rightarrow "$ If the algebra $\Bbb{H}_{\Bbb{Q}%
}\left( -13,q\right) $ is a split algebra, similarly with i), we obtain that for $q\neq 13:$ $\left(\frac{-13,q}{q}\right)=1,$ $\left(\frac{-13,q}{13}\right)=1,$ $\left(\frac{-13,q}{2}\right)=1.$\\
Using the properties of Hilbert symbol and Legendre symbol, similarly with i), we obtain that 
$q\equiv 1; 9; 17; 25; 29; 49$ (mod $52$).\\
If $q=13,$ similarly with i) The algebra  $\Bbb{H}_{\Bbb{Q}%
}\left( -13,13\right) $ is a split algebra.
\begin{corollary}\label{1.} 
Let $q$ be an odd positive prime integer and let $K$ be an algebraic numbers field such that $\left[K:\mathbb{Q}\right]$ is odd. The following statements hold true:\\
i) the algebra $\Bbb{H}_{K
}\left( -3,q\right) $ splits if and only if the algebra $\Bbb{H}_{\Bbb{Q}%
}\left( -3,q\right) $ splits if and only if $q=x^{2}+3y^{2}$ for some $x,y$ $\in$ $\mathbb{Z}$;\\
ii) the algebra $\Bbb{H}_{K
}\left( -5,q\right) $ splits if and only if the algebra $\Bbb{H}_{\Bbb{Q}%
}\left(-5,q\right) $ splits if and only if $q=x^{2}+5y^{2}$ for some $x,y$ $\in$ $\mathbb{Z}$;\\
iii) the algebra $\Bbb{H}_{K
}\left( -7,q\right) $ splits if and only if
the algebra $\Bbb{H}_{\Bbb{Q}%
}\left( -7,q\right) $ splits if and only if $q=x^{2}+7y^{2}$ for some $x,y$ $\in$ $\mathbb{Z}$;\\
iv) the algebra $\Bbb{H}_{K
}\left( -13,q\right) $ splits if and only if
the algebra $\Bbb{H}_{\Bbb{Q}%
}\left(-13,q\right) $ splits if and only if $q=x^{2}+13y^{2}$ for some $x,y$ $\in$ $\mathbb{Z}$.
\end{corollary}
\textbf{Proof.} The proof is immediate using Theorem 1.2., Proposition 2.5 and Remark 1.1.
\vspace{2mm}\\
A question which can appears is the following: when in general a quaternion algebra $\Bbb{H}_{\Bbb{Q}%
}\left( -n,q\right) $ (where $n,q$ $\in$ $\mathbb{N}^{*},$ $q$ is a prime number) is a split algebra?
In [16] we find the following result: "let an odd prime $q$ and $n$$\in$ $\mathbb{Z}$ such that $q-n$ is a square. Then 
$\Bbb{H}_{\Bbb{Q}%
}\left( n,q\right) $ splits if and only if $q\equiv 1$ (mod $4$)"'.\\
In the future, we will study when a quaternion algebra $\Bbb{H}_{\Bbb{Q}%
}\left( n,q\right), $ where $q-n$ is not a square, is a split algebra.
\vspace{2mm}\\
Let $q$ be an odd prime positive integer. Let $K$ be an algebraic number field and $p$ be a prime (finite of infinite) of $K.$ Let $K_{p}$ be the completion of $K$ with the respect to $p$- adic valuation and let $\xi$ be a primitive root of order $q$ of unity such that $\xi$$\in$$K_{p}.$ We consider the symbol algbebra $A=\left(\frac{\alpha,\beta}{K_{p},\xi}\right),$ $\alpha,\beta$$\in$$K^{*}_{p}.$
In the paper [21], we determined some symbol algebras of degree $q=3$ over a local field, which are split algebras. 
\begin{proposition}\label{p2.6.} ([21])
Let $p$ be a prime positive integer, $p\equiv 2$ (\textit{mod}
$3$)
and let the $K_{p}-$ algebra $A=\left( \frac{%
\begin{array}{c}
\\
\alpha ,p^{3l}
\end{array}
}{%
\begin{array}{c}
K_{p},\varepsilon
\end{array}
}\right) ,$ where $\varepsilon$ is a primitive root of order $3$ of the unity, $l\in \Bbb{N}^{*},\alpha \in
K,K=\Bbb{Q}\left( \varepsilon \right) .$ Let $P$ 
be a prime ideal of the ring of integers of the field 
$L=K\left(\sqrt[3]{\alpha}\right),$ lying above
$p.$ Then $p^{3l}$ is a
norm from $L_{P}^{*}$ and the local Artin symbol $\left( \frac{%
\begin{array}{c}
\\
L_{P}\,/\,K_{p}
\end{array}
}{%
\begin{array}{c}
(p^{3l})
\end{array}
}\right) $ is the identity.
\end{proposition}

\begin{proposition}\label{p2.7.} ([21])
Let $p$ be a prime positive integer, $p\equiv 1\,\,$(\textit{%
mod\thinspace \thinspace }$3$)and let $K_{p_{1}}-$%
 algebra $A=\left( \frac{%
\begin{array}{c}
\\
\alpha ,p^{3l}
\end{array}
}{%
\begin{array}{c}
K_{p_{1}},\varepsilon
\end{array}
}\right) ,$ where $\varepsilon$ is a primitive root of order $3$ of the unity, $l\in \Bbb{N}^{*},\alpha \in
K,K=\Bbb{Q}\left(
\varepsilon \right)$ and $p_{1}$ is a prime element in $%
\Bbb{Z}[\varepsilon ],$\textit{\ }$p_{1}\mid p.$ Let
$P$
be a prime ideal in the ring of integers of the field  $L=K\left(\sqrt[3]{\alpha}\right),$ lying above $p_{1}.$ Then $%
p^{3l}\in N_{L_{P}/K_{p_{1}}}\left( L_{P}^{*}\right) $ and
the
local Artin symbol $\left( \frac{%
\begin{array}{c}
\\
L_{P}\,/\,K_{p_{1}}
\end{array}
}{%
\begin{array}{c}
(p^{3l})
\end{array}
}\right) $ is the identity in the Galois group $Gal\left(
L_{P}/K_{p_{1}}\right) .$
\end{proposition}
Now we generalise these results, finding some symbol algebras of degree $q\in\left\{3,5,7,11,13,17,19\right\}$ over local fields, which are split algebras. We obtain the following result:

\begin{proposition}\label{p2.8.}
Let $p$ be a prime positive integer and $q$ $\in$ $\left\{3,5,7,11,13,17,19\right\}.$ Let $\xi$ be a primitive root of order $q$ of the unity and the cyclotomic field $K=\mathbb{Q}\left(\xi\right).$ Let $K_{p_{1}}-$%
 algebra $A=\left( \frac{%
\begin{array}{c}
\\
\alpha ,p^{ql}_{1}
\end{array}
}{%
\begin{array}{c}
K_{p_{1}},\xi
\end{array}
}\right) ,$ where $l\in \Bbb{N}^{*},\alpha \in
K$ and $p_{1}$ is a prime element in $%
\Bbb{Z}[\xi ],$\textit{\ }$p_{1}\mid p.$ Let
$P$
be a prime ideal in the ring of integers of the field $L=K\left(\sqrt[q]{\alpha}\right),$ lying above $p_{1}.$ Then $%
p^{ql}_{1}\in N_{L_{P}/K_{p_{1}}}\left( L_{P}^{*}\right) $ 
\end{proposition}
\textbf{Proof.} We denote with $\mathcal{O}_L$ the ring of integers of the field $L=K\left(\sqrt[q]{\alpha}\right).$\\
\textbf{Case 1.} If $<\overline{p}>=\left(\mathbb{Z}^{*}_{q},\cdot\right),$ from Proposition 1.3, we obtain that $p$ is a prime in the ring $\mathbb{Z}\left[\xi\right],$ therefore $p_{1}=p$ and the $q$- power character $\left(\frac{\alpha}{p\mathbb{Z}\left[\xi\right]}\right)_{q}=1.$ Applying Proposition 1.4 it results that $p$ is totally split in $\mathcal{O}_L:$ $p\mathcal{O}_L=P_{1}P_{2}\cdot...\cdot P_{q},$ where $P_{i}$$\in$Spec($\mathcal{O}_L$), $i=\overline{1,q}.$ \\
If we denote with $g$ the number of decomposition of the ideal 
$p\mathcal{O}_L,$ with $e_{i}$ the ramification index of $p$ at $P_{i}$ and with $f_{i}=$$[\mathcal{O}_L/P_{i} : \mathcal{O}_K/p\mathcal{O}_K]$ the residual degree of p ($i=\overline{1,q}$), since the fields extension $K\subset L$ is a Galois extension we have $e_{1}=e_{2}=...=e_{q}=e,$ $f_{1}=f_{2}=...=f_{q}=f$ and $efg=\left[L:K\right]=q.$ But $g=q,$ therefore $e=f=1.$ It is known that $\left[L_{P}:K_{p}\right]=ef,$ then $L_{P}=K_{p},$ for each $P$$\in$Spec($\mathcal{O}_L$), $P|p\mathcal{O}_L.$ We obtain that $p$ is the norm of itself in the trivial extension $K_{p}\subseteq L_{P}.$\\
\textbf{Case 2.} If $<\overline{p}>\neq\left(\mathbb{Z}^{*}_{q},\cdot\right),$ applying Proposition 1.3, we obtain that the number of decomposition of the ideal $p\mathbb{Z}\left[\xi\right]$ in the product of prime ideals in the ring $\mathbb{Z}\left[\xi\right]$ is $g^{'}=\frac{\varphi\left(q\right)}{f^{'}}=\frac{q-1}{f^{'}},$ where $f^{'}=ord_{\left(\mathbb{Z}^{*}_{q},\cdot\right)}\overline{p}.$ Using Proposition 1.2, it results that there are $p_{1}, p_{2},...,p_{g^{'}},$ prime elements in $\mathbb{Z}\left[\xi\right]$ such that $p\mathbb{Z}\left[\xi\right]=p_{1}\mathbb{Z}\left[\xi\right]\cdot p_{2}\mathbb{Z}\left[\xi\right]\cdot...\cdot p_{g^{'}}\mathbb{Z}\left[\xi\right].$\\
\textbf{Subcase 2 a).} If the $q-$ power character $\left(\frac{\alpha}{p_{1}\mathbb{Z}\left[\xi\right]}\right)_{q}=1,$ similarly with the case 1, we obtain that $p_{1}$ is a norm of itself in the trivial extension $K_{p_{1}}\subseteq L_{P},$ where $P$ is a prime ideal of $\mathcal{O}_L$ lying above $p_{1}.$\\
\textbf{Subcase 2 b).} If the $q-$ power character $\left(\frac{\alpha}{p_{1}\mathbb{Z}\left[\xi\right]}\right)_{q}$ is a root of order $q$ of unity different fromn $1,$ applying Proposition 1.4 we obtain that $p_{1}\mathcal{O}_L$$\in$ Spec($\mathcal{O}_L$). Knowing that $efg=\left[L:K\right]=q$ and $g=e=1,$ it results that the residual degree is $f=q,$ therefore $f|ql.$ Applying Theorem 1.3, we obtain that $p^{ql}_{1}$$\in$$N_{L_{P}/K_{p_{1}}}\left(L^{*}_{P}\right).$
\vspace{0.1in}

\textbf{Acknowledgements.} The author thanks Professor Victor Alexandru for helpful discussions on this topic.

\vspace{2mm} \noindent \footnotesize
\begin{minipage}[b]{10cm}
Diana SAVIN, \\
Department of Mathematics and Computer Science, \\ 
Ovidius University of Constanta, \\ 
Constanta 900527, Bd. Mamaia no.124, Rom\^{a}nia \\
Email: savin.diana@univ-ovidius.ro\end{minipage}

\end{document}